\documentclass[12pt]{article}
\usepackage{amssymb}
\def\bbba{\overline{\mathbb Q}}
\def\bbbq{{\mathbb Q}}
\def\bbbc{{\mathbb C}}

\def\qed{\hfill{\bf qed}}
\def\trace{{\rm Trace}}
\newtheorem{theorem}[subsection]{Theorem}

\newtheorem{lemma}[subsection]{Lemma}

\newtheorem{corollary}[subsection]{Corollary}
\newtheorem{proposition}[subsection]{Proposition}

\begin{document}
\title{A refined version of the Siegel-Shidlovskii theorem}
\author{F.Beukers}
\maketitle

\abstract{Using Y.Andr\'e's result on differential equations
satisfied by $E$-functions, we derive an improved version of
the Siegel-Shidlovskii theorem. It gives a complete characterisation
of algebraic relations over the algebraic numbers between
values of $E$-functions at any non-zero algebraic point.}

\section{Introduction}
In this paper we consider $E$-functions. An entire function
$f(z)$ is called an $E$-function if it has a 
powerseries expansion of the form
$$f(z)=\sum_{k=0}^{\infty} {a_k\over k!} z^k$$
where 
\begin{enumerate}
\item $a_k\in\bbba$ for all $k$.
\item $h(a_0,a_1,\ldots,a_k)=O(k)$ for all $k$ where $h$
denotes the log of the absolute height.
\item $f$ satisfies a linear differential equation $Ly=0$
with coefficients in $\bbba[z]$.
\end{enumerate}

The linear differential equation $Ly=0$ of minimal order which is
satisfied by $f$ is called the {\it minimal differential equation}
of $f$.

Furthermore, in all of our consideration we take a fixed
embedding $\bbba\to \bbbc$.

Siegel first introduced $E$-functions around 1929 in his
work on transcendence of values of Bessel-functions and
related functions. Actually, Siegel's definition was slightly
more general in that condition (3) reads 
$h(a_0,a_1,\ldots,a_k)=o(k\log k)$. But until now no $E$-functions
in Siegel's original definition are known which fail to satisfy
condition (2) above. Around 1955 Shidlovski managed to remove
Siegel's technical normality conditions and we now have the
following theorem (see [Sh, Chapter 4.4],[FN, Theorem 5.23]). 

\begin{theorem}[Siegel-Shidlovskii, 1956]\label{siegelshidlovskii}
Let $f_1,\ldots,f_n$ be
a set of $E$-functions which satisfy the system of first
order equations
$${d\over dz}\pmatrix{y_1\cr \vdots\cr y_n\cr}=A\pmatrix{y_1\cr
\vdots\cr y_n\cr}$$
where $A$ is an $n\times n$-matrix with entries in $\bbba(z)$.
Denote the common denominator of the entries of $A$ by $T(z)$.
Then, for any $\xi\in\bbba$ such that $\xi T(\xi)\ne0$ we have
$${\rm degtr}_{\bbba}(f_1(\xi),\ldots,f_n(\xi))
={\rm degtr}_{\bbba(z)}(f_1(z),\ldots,f_n(z)).$$
\end{theorem}

In [B1] Daniel Bertrand gives an alternative proof of the
Siegel-Shidlovskii theorem using Laurent's determinants.

Using the Siegel-Shidlovskii Theorem it is possible to prove
the following theorem.

Of course the Siegel-Shidlovskii theorem suggests strongly that
all relations between values of $E$-functions at algebraic
points arise by specialisation of polynomial relations over
$\bbba(z)$. Using the techniques of the Siegel-Shidlovskii
techniques this turns out to be true up to a finite exceptional
set of algebraic points.

\begin{theorem}[Nesterenko-Shidlovskii, 1996]\label{specialisation} 
There exists a finite set $S$ such that for
all $\xi\in\bbba,\xi\not\in S$ the following holds. 
For any homogeneous polynomial relation $P(f_1(\xi),\ldots,f_n(\xi))=0$
with $P\in\bbba[X_1,\ldots,X_n]$ there exists 
$Q\in\bbba[z,X_1,\ldots,X_n]$, homogeneous in $X_i$,
such that $Q(z,f_1(z),\ldots,f_n(z))
\equiv0$ and $P(X_1,\ldots,X_n)=Q(\xi,X_1,\ldots,X_n)$.
\end{theorem}

In the statement of the Theorem one can drop the word `homogeneous'
if one wants, simply by considering the set of $E$-functions
$1,f_1(z),\ldots,f_n(z)$ instead. 
Loosely speaking, for almost all $\xi\in\bbba$, polynomial relations
between the values of $f_i$ at $z=\xi$ arise by specialisation of
polynomial relations between the $f_i(z)$ over $\bbba(z)$.

In [NS] it is also remarked that the exceptional set $S$ can be
computed in principle. Although Theorem
\ref{specialisation} is not stated explicitly in [NS],
it is immediate from Theorem 1 and Lemmas 1,2 in [NS]. 

Around 1997 Y.Andr\'e (see [A1] and Theorem \ref{annihilate} below)
discovered that the nature of differential
equations satisfied by $E$-functions is very simple. Their
only non-trivial singularities are at $0,\infty$. Even more
astounding is that this observation allowed Andr\'e to prove transcendence
statements, as illustrated in Theorem \ref{vanishat1}. In particular
Andr\'e managed to give a completely new proof of the Siegel-Shidlovskii
Theorem using his discovery. In order to achieve this, a defect relation
for linear equations with irregular singularities had to be invoked.
For a survey one can consult [A2] or, more detailed, [B2].

However, it turns out that even more is possible.
Theorem \ref{annihilate} allows us to prove
the following Theorem.

\begin{theorem}\label{specialisation2}
Theorem \ref{specialisation} holds for any $\xi\in\bbba$ with
$\xi T(\xi)\ne0$.
\end{theorem}

The proof of this Theorem
will be given in section 3, after the necessary
preparations. In particular we will use some very basic facts about
differential galois groups of systems of differential equations.
All that we require is contained in Section 1.4 of the book [PS].

One particular consequence of Theorem \ref{specialisation2} is the
solution of Conjecture A in [NS].

\begin{corollary} Let assumption be as in Theorem \ref{siegelshidlovskii}.
Suppose that $f_1(z),\ldots,f_n(z)$ are linearly independent over
$\bbba(z)$. Then for any $\xi\in\bbba$, with $\xi T(\xi)\ne0$,
the numbers $f_1(\xi),\ldots,f_n(\xi)$ are $\bbba$-linear independent.
\end{corollary}

A question that remains is about the nature of relations between
values of $E$-functions at singular points $\ne0$. The best known
example is $f(z)=(z-1)e^z$. Its differential equation has a singularity
at $z=1$ and it vanishes at $z=1$, even though $f(z)$ is transcendental
over $\bbba(z)$. Of course the vanishing of $f(z)$ at $z=1$ arises
in a trivial way and one would probably agree that it is better
to look at $e^z$ itself. It turns out that all relations between
values of $E$-functions at singularities $\ne0$ arise in a similar trivial
fashion. This is a consequence of the following Theorem.

\begin{theorem}\label{essentialE}
Let $f_1,\ldots,f_n$ be as above and suppose they are
$\bbba(z)$-linear independent. Then there exist $E$-functions
$e_1(z),\ldots,e_n(z)$ and an $n\times n$-matrix $M$ with entries
in $\bbba[z]$ such that
$$\pmatrix{f_1(z)\cr \vdots \cr f_n(z)}=M
\pmatrix{e_1(z)\cr \vdots \cr e_n(z)}$$
and where $(e_1(z),\ldots,e_n(z))$ is vector solution of a system
of $n$ homogeneous first order equations with coefficients
in $\bbba[z,1/z]$.
\end{theorem}

{\bf Acknowledgement} At this point I would like to express
my gratitude to Daniel Bertrand who critically read and
commented on a first draft of this paper. His remarks 
were invaluable to me.

\section{Andr\'e's Theorem and first consequences}

Everything we deduce in this paper hinges on the following
beautiful Theorem plus Corollary by Yves Andr\'e.

\begin{theorem}[Y.Andr\'e]\label{annihilate}
Let $f$ be an $E$-function and let $Ly=0$ be its minimal differential
equation. Then at every point $z\ne 0,\infty$ the equation has a
basis of holomorphic solutions.
\end{theorem}

All results that follow now, depend on a limited version of Theorem 
\ref{annihilate} where the $E$-function has rational coefficients.
Although the following theorem occurs in [A1] we like to give a proof
of it to make this paper selfcontained to the extent of only
accepting Theorem \ref{annihilate}.

\begin{corollary}[Y.Andr\'e]\label{vanishat1}
Let $f$ be an $E$-function with rational coefficients and let $Ly=0$ be
its minimal differential equation. Suppose $f(1)=0$. Then
$z=1$ is an apparent singularity of $Ly=0$.
\end{corollary}

{\bf Proof}
Suppose 
$$f(z)=\sum_{n=0}^{\infty}{a_n\over n!}z^n.$$
Let $g(z)=f(z)/(1-z)$. Note that $g(z)$ is also holomorphic
in $\bbbc$. Moreover, $g(z)$ is again an $E$-function. 
Write
$$g(z)=\sum_{n=0}^{\infty}{b_n\over n!}z^n$$
where
$${b_n\over n!}=\sum_{k=0}^n{a_k\over k!}.$$
Since $f(1)=0$ we see that 
$${b_n\over n!}=-\sum_{k=n+1}^{\infty}{a_k\over k!}.$$
Since $f$ is an $E$-function there exist $B,C>0$ such
that $|a_k|\le B\cdot C^k$. Hence
\begin{eqnarray*}
|b_n|&\le& Bn!\left|\sum_{k=n}^{\infty}{C^k\over k!}\right|\\
&\le& Bn!{C^n\over n!}\left(1+{C\over 1!}+{C^2\over 2!}+\cdots
\right)\\
&\le&Be^C\cdot C^n
\end{eqnarray*}
Furthermore, the common denominator of $b_0,\ldots,b_n$ is bounded
above by the common denominator of $a_0,a_1,\ldots,a_n$, hence
bounded by $B_1\cdot C_1^n$ for some $B_1,C_1>0$. This shows that
$f(z)/(z-1)$ is an $E$-function. 
The minimal differential operator which annihilates $g(z)$ is
simply $L\circ(z-1)$. From Andr\'e's theorem \ref{annihilate}
it follows that
the kernel of $L\circ(z-1)$ around $z=1$ is
spanned by holomorphic functions. Hence the kernel of $L$ is
spanned by holomorphic solutions times $z-1$. In other words,
all solutions of $Ly=0$ vanish at $z=1$ and therefore
$z=1$ is an apparent singularity. 

\qed
\medskip

\begin{lemma}\label{connectedcomponent}
Let $f$ be an $E$-function with minimal differential equation 
$Ly=0$ of order $n$.
Let $G$ be its differential galois group and let $G^o$
be the connected component of the identity in $G$. Let $V$ be
the vectorspace spanned by all images of $f(z)$ under $G^o$.
Then $V$ is the complete solution space of $Ly=0$.
\end{lemma}

{\bf Proof}
The fixed field of $G^o$ is an algebraic Galois extension $K$ of $\bbba(z)$
with galois group $G/G^o$.
Suppose that $V$ has dimension $m$. Then $f$ satisfies a linear
differential equation with coefficients in $K$ of order $m$.
In particular we have a relation  
\begin{equation}\label{relation}
f^{(m)}+p_{m-1}(z)f^{(m-1)}+\cdots+p_1(z)f'+p_0(z)f=0
\end{equation}
for some $p_i\in K$. We subject this relation to analytic continuation.
Since $f$ is an entire function, it has trivial monodromy. By
choosing suitable paths we obtain the conjugate relations
$$f^{(m)}+\sigma(p_{m-1})f^{(m-1)}+\cdots+\sigma(p_1)f'+\sigma(p_0)f=0$$
for all $\sigma\in G/G^o$.
Taking the sum over all these relations gives us a non-trivial differential
equation for $f$ of order $m$ over $\bbba(z)$.
From the minimality of $Ly=0$ we now conclude that $m=n$, i.e.
the dimension of $V$ is $n$.

\qed
\medskip

Actually it follows from Theorem \ref{annihilate} that the fixed 
field of $G^o$ is of the form $K=\bbba(z^{1/r})$ for some positive
integer $r$. But we don't need that in our proof.

The following Lemma is straightforward consequence of
general theory of algebraic groups.

\begin{lemma}\label{productconnect}
Let $G_1,\ldots,G_r$ be linear algebraic groups and denote
by $G_i^o$ their components of the identity. Let 
$H\subset G_1\times G_2\times\cdots\times G_r$ be an algebraic
subgroup such that the natural projection $\pi_i:H\to G_i$
is surjective for every $i$. Let $H^o$ be the connected component
of the identity in $H$. Then the natural projections $\pi_i:H^o\to
G_i^o$ are surjective.
\end{lemma}

Now we prove a generalistaion of Andr\'e's Corollary \ref{vanishat1}
to general non-zero algebraic points.

\begin{theorem}\label{vanishatalgebraic}
Let $f$ be an $E$-function with minimal differential equation
$Ly=0$ of order $n$. Suppose that $\xi\in\bbba^*$ and $f(\xi)=0$.
Then all solutions of $Ly$ vanish at $z=\xi$. In particular,
$Ly=0$ has an apparent singularity at $z=\xi$.
\end{theorem}

{\bf Proof} By replacing $f(z)$ by $f(\xi z)$ if necessary,
we can assume that $f$ vanishes at $z=1$.
Let $f^{\sigma_1}(z),\ldots,f^{\sigma_r}(z)$ be the 
$Gal(\bbba/\bbbq)$-conjugates
of $f(z)$ where we take $f^{\sigma_1}(z)=f(z)$. 
Let $L^{\sigma_i}y=0$ be the $\sigma_i$-conjugate of $Ly=0$.
Note that this is the minimal differential equation satisfied
by $f^{\sigma_i}(z)$. Let $G_i$ be the differential galois
group and $G_i^o$ the connected component of the identity.
By Lemma \ref{connectedcomponent} the images of $f^{\sigma_i}(z)$
under $G_i^o$ span the
complete solution space of $L^{\sigma_i}y=0$. 

The product
$F(z)=\prod_{i=1}^r f^{\sigma_i}(z)$ is an $E$-function having
rational coefficients. Let ${\cal L}y=0$ be its minimal differential
equation. Furthermore, $F(1)=0$. Hence, from Andr\'e's Theorem
\ref{vanishat1} it follows that all solutions of ${\cal L}y=0$ vanish at $z=1$.

Let $H$ be the differential galois group of the differential compositum of
the Picard-Vessiot extensions corresponding to $L^{\sigma_i}y=0$.
Note that the image of $F(z)$ under any $h\in H$ is again a solution
of ${\cal L}y=0$. In particular this image also vanishes at $z=1$.

Furthermore, $H$
is an algebraic subgroup of $G_1\times G_2\times\cdots\times G_r$
such that the natural projections $\pi_i:H\to G_i$ are surjective.
Let $H^o$ be the connected component of the identity of $H$. Then,
by Lemma \ref{productconnect},
the projections $\pi_i:H^o\to G_i^o$ are surjective. 

Let $V_i$ be the solution space of local solutions at $z=1$ of
$L^{\sigma_i}y=0$. Let $W_i$ be the linear subspace of solutions
vanishing at $z=1$. The group $H^o$ acts linearly on each space
$V_i$. Let $v_i\in V_i$ be the vector corresponding to the
solution $f^{\sigma_i}(z)$. Define $H_i=\{h\in H^o| \pi_i(h)v_i\in W_i\}$.
Then $H_i$ is a Zariski closed subset of $H^o$. Furthermore, because
all solutions of ${\cal L}y=0$ vanish at $z=1$, we have
that $H^o=\cup_{i=1}^r H_i$. Since $H^o$ is connected this implies that
$H_i=H^o$ for at least one $i$. Hence $\pi_i(H_i)=\pi_i(H^o)=G_i^o$
and we see that $gv_i\in W_i$ for all $g\in G_i^o$. We conclude that
$W_i=V_i$. In other words, all local solutions of $L^{\sigma_i}y=0$
around $z=1$ vanish in $z=1$. By conjugation we now see that
the same is true for $Ly=0$. 

\qed
\section{Independence results}
We now consider a set of $E$-functions $f_1,\ldots,f_n$ which
satisfy a system of homogeneous first order equations
$$y'=Ay$$
where $y$ is a vector of unknown functions $(y_1,\ldots,y_n)^t$
and $A$ an $n\times n$-matrix with entries in $\bbba(z)$. The
common denominator of these entries is denoted by $T(z)$.

\begin{lemma}\label{independentrelations}
Let us assume that the $\bbba(z)$-rank of $f_1,\ldots,f_n$ is $m$.
Then there is a $\bbba[z]$-basis of relations
\begin{equation}\label{linrelations}
C_{i,1}(z)f_1(z)+C_{i,2}(z)f_2(z)+\cdots+C_{i,n}(z)f_n(z)\equiv 0
\quad i=1,2,\ldots, n-m
\end{equation}
such that for any $\xi\in\bbba$ the matrix
$$\pmatrix{C_{11}(\xi) & C_{12}(\xi) & \ldots & C_{1n}(\xi)\cr
\vdots & \vdots & & \vdots \cr
C_{n-m,1}(\xi) & C_{n-m,2}(\xi) & \ldots & C_{n-m,n}(\xi)\cr}$$
has rank precisely $n-m$. 
\end{lemma}

{\bf Proof} The $\bbba(z)$-dimension of all relations is $n-m$.
Choose an indepent set of $n-m$ relations of the form
(\ref{linrelations}) (without the extra specialisation condition).

Denote the greatest common divisor of the determinants of 
all $(n-m)\times(n-m)$ submatrices of $(C_{ij}(z))$ by $D(z)$.
Suppose that $D(\xi)=0$ for some $\xi$. Then the matrix $(C_{ij}(\xi))$
has linearly dependent rows. By taking $\bbba$-linear relations 
between the rows, if necessary, we can assume that $C_{1j}(\xi)=0$
for $j=1,\ldots,n$. Hence all $C_{1j}(z)$ are divisible by $z-\xi$
and the polynomials $C_{1j}(z)/(z-\xi)$ are the coefficients of another
$\bbba(z)$-linear relation. Replace the first relation by this new relation.
The new greatest divisor of all $(n-m)\times(n-m)$-determinants
is now $D(z)/(z-\xi)$. By repeating this argument we can find
an independent set of $n-m$ relations of the form (\ref{linrelations})
whose associated $D(z)$ is a non-zero constant. 

But now it is not hard to see that (\ref{linrelations}) is a 
$\bbba[z]$-basis of all $\bbba[z]$-relations. Furthermore, $D(\xi)\ne0$
for all $\xi$ (because $D(z)$ is constant), so all specialisations
have maximal rank.

\qed

\begin{theorem}\label{relations}
Let $f_1,\ldots,f_n$ be a vector solution of the system
$$y'=Ay$$
consisting of $E$-functions. Let $T(z)$ be the common denominator
of the entries in $A$. 
Then, for any $\xi\in\bbba$, $\xi T(\xi)\ne0$, any $\bbba$-linear
relation between $f_1(\xi),\ldots,f_n(\xi)$ arises by specialisation
of a $\bbba(z)$-linear relation.
\end{theorem}

{\bf Proof}
Suppose there exists a $\bbba$-linear relation
$$\alpha_1f_1(\xi)+\alpha_2f_2(\xi)+\cdots+\alpha_nf_n(\xi)=0$$
which does not come from specialisation of a $\bbba(z)$-linear
relation at $z=\xi$. Consider the function
$$F(z)=A_1(z)f_1(z)+A_2(z)f_2(z)+\cdots+A_n(z)f_n(z)$$
where $A_i(z)\in \bbba[z]$ to be specified later. Let $Ly=0$
be the minimal differential equation satisfied by $F$. Suppose
that the $\bbba(z)$-rank of $f_1,\ldots,f_n$ is $m$. It will
turn out that the order of $Ly=0$ is at most $m$. 

We now show how to choose $A_1(z),\ldots,A_n(z)$  such that
\begin{enumerate}
\item $A_i(\xi)=\alpha_i$ for $i=1,2,\ldots,n$
\item The order of $Ly=0$ is $m$.
\item $\xi$ is a regular point of $Ly=0$.
\end{enumerate}

By using the system $y'=Ay$ recursively we can find
$A_i^j(z)\in\bbba[z]$ such that 
$$F^{(j)}(z)=\sum_{i=1}^n A_i^j(z)f_i(z).$$
In addition we fix a $\bbba(z)$-basis of linear relations
$$C_{i,1}(z)f_1(z)+\cdots+C_{i,n}f_n(z)\equiv0\quad i=1,\ldots,n-m$$
with polynomial coefficients $C_{ij}(z)$
such that the $(n-m)\times n$-matrix of values $C_{ij}(\xi)$
has maximal rank $n-m$. This is possible in view
of Lemma \ref{independentrelations}. Consider the $(n+1)\times n$-matrix
$${\cal M}=\pmatrix{C_{11}(z) & \ldots & C_{1n}(z)\cr
\vdots & & \vdots\cr
C_{n-m,1}(z) & \ldots & C_{n-m,n}(z)\cr
A_1(z) & \ldots & A_n(z)\cr
\vdots & & \vdots\cr
A_1^m(z) & \ldots & A_n^m(z)\cr}.$$
We denote the submatrix obtained from ${\cal M}$ by deleting
the row with $A_i^j\ (i=1,\ldots,n)$ by ${\cal M}_j$. 
There exists a $\bbba(z)$-linear relation between the rows
of ${\cal M}$ which explains why the minimal equation $Ly=0$ of
$F$ satisfies a differential
equation of order $\le m$. This equation has precisely order
$m$ if and only if the submatrix ${\cal M}_m$ has rank $m$. In that
case the minimal differential equation for $F$ is given by
$$\Delta_m F^{(m)}+\ldots+\Delta_1 F'+\Delta_0 F=0$$
where $\Delta_j=(-1)^j\det({\cal M}_j)$.

By induction it is not hard to show that $A_i^0(z)=A_i(z)$ and
$$A_i^j(z)=A_i^{(j)}+P_{ij}(A_1,\ldots,A_n,\ldots,A_1^{(j-1)},
\ldots,A_n^{(j-1)})$$
where 
$$P_{ij}\in\bbba[z,1/T(z)][X_{10},\ldots,X_{n0},\ldots,X_{1,j-1},
\ldots,X_{n,j-1}]$$
are linear forms with coefficients in $\bbba[z,1/T(z)]$.
We can now choose the $A_i(z)$ and their derivatives in such a
way that $\det({\cal M}_m)$ does not vanish in the point $\xi$.
The choice of $A_i(\xi)$ is fixed by taking $A_i(\xi)=\alpha_i$.
Since the relation $\sum_{i=1}^n \alpha_if_i(\xi)=0$ does not come
from specialisation, the rows of values $(C_{i1}(\xi),\ldots,C_{in}(\xi))$
for $i=1,\ldots,n-m$ and $(\alpha_1,\ldots,\alpha_n)$ have maximal
rank $n-m+1$. We can now choose the derivates $A_i^{(j)}$ recursively with
respect to $j$ such that $\det({\cal M}_m)(\xi)\ne0$. With this choice
we note that conditions (i),(ii),(iii) are satisfied. 

On the other hand, $F(\xi)=0$, so it follows from 
Theorem \ref{vanishatalgebraic}
that $\xi$ is a singularity
of $Ly=0$. This contradicts condition (iii).

\qed

{\bf Proof of Theorem \ref{specialisation2}}.
Consider the vector of $E$-functions given by the monomials
${\bf f}(z)^{\bf i}:=f_1(z)^{i_1}\cdots f_n(z)^{i_n}$, $i_1+\ldots+i_n=N$ of
degree $N$ in $f_1(z),\ldots,f_n(z)$. This vector again satisfies a
system of linear first order equations with singularities in
the set $T(z)=0$. So we now apply Theorem \ref{relations} to
the set of $E$-functions ${\bf f}(z)^{\bf i}$. The relation
$P(f_1(\xi),\ldots,f_n(\xi))$ is now a $\bbba$-linear relation
between the values ${\bf f}(\xi)^{\bf i}$. Hence, by Theorem 
\ref{relations}, there is a $\bbba[z]$-linear relation between the
${\bf f}(z)^{\bf i}$ which specialises to the linear relation
between the values at $z=\xi$. This proves our Theorem.

\qed

\section{Removal of non-zero singularities}
In this section we prove Theorem \ref{essentialE}.
For this we require the following Proposition.

\begin{proposition}\label{quotientisE}
Let $f$ be an $E$-function and $\xi\in\bbbq^*$ such that $f(\xi)=0$.
Then $f(z)/(z-\xi)$ is again an $E$-function.
\end{proposition}

{\bf Proof} By replacing $f(z)$ by $f(\xi z)$ if necessary,
we can restrict our attention to $\xi=1$.
Write down a basis of local solutions of $Ly=0$ around $z=1$.
Since $f$ vanishes at $z=1$, Theorem \ref{vanishatalgebraic}
implies that all solutions of $Ly=0$ vanish at $z=1$. But then,
by conjugation, this holds for the solutions around $z=1$ of
the $Gal(\bbba/\bbbq)$-conjugates $L^{\sigma}y=0$ as well.
In particular, the conjugate $E$-function $f^{\sigma}(z)$ 
vanishes at $z=1$ for every $\sigma\in Gal(\bbba/\bbbq)$.
Taking up the notations of the proof of Theorem \ref{vanishat1}
we now see that
$${b_n^{\sigma}\over n!}=-\sum_{k=n+1}^{\infty}{a_k^{\sigma}\over k!}$$
for every $\sigma$. We can now bound $|b_n^{\sigma}|$ exponentially
in $n$ for every $\sigma$. Since the coefficients of an $E$-function
lie in a finite extension of $\bbbq$, only finitely many conjugates
are involved. So we get our desired bound $h(b_0,\ldots,b_n)=O(n)$.

\qed
\medskip

{\bf Proof of Theorem \ref{essentialE}}
Denote the column vector $(f_1(z),\ldots,f_n(z))^t$ by 
${\bf f}(z)$. Let 
$${\bf y}'(z)=A(z){\bf y}(z)$$
be the system of equations satisfied by ${\bf f}$
and let $G$ be its differential Galois group.
Because the $f_i(z)$ are $\bbba(z)$-linear independent, the
images of ${\bf f}$ under $G$ span the complete solution set
of ${\bf y}'=A{\bf y}$. So the images under $G$ give us
a fundamental solution set ${\cal F}$ of our system.
We assume that the first column is ${\bf f}(z)$ itself.
Since the $f_i(z)$
are $E$-functions, it follows from Theorem \ref{annihilate}
that the entries of ${\cal F}$
are holomorphic at every point $\ne0$. Consequently, the determinant
$W(z)=\det({\cal F})$ is holomorphic outside $0$. Since $W(z)$
satisfies $W'(z)=\trace(A)W(z)$, we see that $W(\alpha)=0$ implies 
that $\alpha$ is a singularity of our system. In particular, $\alpha\in\bbba$.
Let $k$ be the highest order with which $\alpha$ occurs as pole in $A$.
Write $\tilde{A}(z)=(z-\alpha)^kA(z)$. Then it follows from specialisation
at $z=\alpha$ of $(z-\alpha)^k{\bf f}'(z)=\tilde{A}(z){\bf f}(z)$ that
there is a non-trivial vanishing relation between the components
of ${\bf f}(\alpha)$. By choosing a suitable $M\in GL(n,\bbba)$
we can see to it that $M{\bf f}(z)$ is a vector of $E$-functions, of
which the first component vanishes at $\alpha$. But then, by
Theorem \ref{vanishatalgebraic}, the whole first row of $M{\cal F}(z)$
vanishes at $z=\alpha$. Hence we can write 
$M{\cal F}(z)=D{\cal F}_1$ where
$$D=\pmatrix{z-\alpha & 0 & \ldots & 0\cr 0 & 1 & \ldots &0\cr
\vdots & \vdots & &\vdots\cr
0 & 0 & \ldots & 1\cr}$$
and ${\cal F}_1$ has entries holomorphic around $z=\alpha$. Thanks
to Proposition \ref{quotientisE}, the entries of the first column
in ${\cal F}_1$ are again $E$-functions. Moreover, ${\cal F}_1$ satisfies
the new system of equations
$${\cal F}_1'=(D^{-1}MAM^{-1}D+D^{-1}D'){\cal F}_1.$$
Notice that the order of vanishing of $W_1(z)=\det({\cal F}_1)$ at
$z=\alpha$ is $1$ lower than the vanishing order of $W(z)$.
We repeat our argument when $W_1(\alpha)=0$. By
using this reduction procedure to all zeros of $W(z)$ we end up
with an $n\times n$-matrix $B$, with entries in $\bbba[z]$,
and an $n\times n$-matrix of holomorphic functions ${\cal E}$
such that ${\cal F}=B{\cal E}$, the first column of ${\cal E}$
consists of $E$-functions and $\det({\cal E})$ is nowhere
vanishing in $\bbbc^*$.
As a result we have ${\cal E}'(z)=A_E(z){\cal E}(z)$ where $A_E(z)$
is an $n\times n$-matrix with entries in $\bbba[z,1/z]$.

\qed 
\section{References}
\begin{itemize}
\item[[A1]] Y.Andr\'e, S\'eries Gevrey de type arithm\'etique I,II,
Annals of Math 151(2000), 705-740, 741-756.
\item[[A2]] Y.Andr\'e, Arithmetic Gevrey series and transcendence,
a survey, J.Th\'eorie des Nombres de Bordeaux 15(2003), 1-10.
\item[[B1]] D.Bertrand, Le th\'eor\`eme de Siegel-Shidlovskii r\'evisit\'e, 
pr\'epublication de l'institut de maths, Jussieu 370 (2004), 
\item[[B2]] D.Bertrand, On Andr\'e's proof of the Siegel-Shidlovskii
theorem. Colloque Franco-Japonnais: Th\'eorie des nombres
transcendants (Tokyo 1998), 51-63, Sem. Math.Sci 27, 
Keio Univ,. Yokohama 1999.
\item[[FN]] N.I.Fel'dman, Yu.V.Nesterenko, {\it Transcendental Numbers}
Encyclopedia of Mathematical Sciences (eds A.N.Parshin, I.R.Shafarevich),
Vol 44 (Number Theory IV, translated by N.Koblitz), Springer Verlag 1998.
\item[[NS]] Yu.V.Nesterenko, A.B.Shidlovskii, On the linear independence
of values of $E$-functions, Mat.Sb 187(1996), 93-108, translated
in Sb.Math 187 (1996), 1197-1211.
\item[[PS]] M. dan der Put, M.Singer, {\it Galois Theory of Differential
Equations}, Vol 328, Grundlehren, Springer Verlag 2003.
\item[[Sh]] A.B.Shidlovski, {\it Transcendental Numbers}, W.de Gruyter 
Studies in Mathematics 12, 1989.
\end{itemize}
\end{document}